\documentclass[10pt,reqno]{amsart}   
\usepackage{amssymb,amscd,latexsym}   
\usepackage{amsmath}
\usepackage{amsthm,times}
\usepackage[all]{xy}
\usepackage{epsfig,graphicx}
\usepackage{graphicx,color}
\usepackage[utf8]{inputenc}
\usepackage{float}
\usepackage{tikz}
\usepackage{enumitem}

\usepackage{tikz-cd}

\newcommand{\llar}{-\kern-5pt-\kern-5pt\longrightarrow}

\usepackage{cases}

\usepackage{multirow}
\usepackage{makecell} 

\newtheorem{Theorem}{Theorem} 
\newtheorem{Lemma}[Theorem]{Lemma}
\newtheorem{cor}[Theorem]{Corollary}
\newtheorem{prop}[Theorem]{Proposition}
\newtheorem{Remark}[Theorem]{Remark}

\DeclareMathOperator{\Hom}{Hom}

\DeclareMathOperator{\Image}{Im}
\DeclareMathOperator{\coker}{coker}
\DeclareMathOperator{\Ker}{Ker}

\DeclareMathOperator{\ext}{Ext}
\DeclareMathOperator{\gl}{GL}

\DeclareMathOperator{\rank}{rank}

\DeclareMathOperator{\E}{\rm E}

\DeclareMathOperator{\EXT}{\mathcal{E}xt}

\newcommand{\op}[1]{{\mathcal O}_{\mathbb{P}^{#1}}}

\def\EE{{\mathcal{E}}}

\def\ker{{\rm ker}\,}

\def\restr{{\kern-1pt\restriction\kern-1pt}}

\def\X{{\mathcal X}}

\def\E{{\mathcal E}}

\newcommand{\simeqd}{\mathrel{\rotatebox[origin=c]{-90}{$\simeq$}}}


\subjclass[2010]{}

\begin{document}

\title[Classification of monads and moduli components] {Classification of monads and moduli components of stable rank 2 bundles with odd determinant and $c_2=10$}

\author{Aislan Leal Fontes}
\address{Departamento de Matem\'atica, UFS - Campus Itabaiana. Av. Vereador Ol\'impio Grande s/n, 49506-036 Itabaiana/SE, Brazil}
\email{aislan@ufs.br}
\author{Marcos Jardim}
\address{IMECC - UNICAMP \\ Departamento de Matem\'atica \\
Rua S\'ergio Buarque de Holanda, 651\\ 13083-970 Campinas-SP, Brazil}
\email{jardim@unicamp.br}

\maketitle

\begin{abstract}
In this paper, we provide a complete classification of the positive minimal 
monads whose cohomology is a stable rank 2 bundle on $\mathbb{P}^3$ with Chern classes $c_1=-1, c_2=10$ and we prove the existence of a new irreducible component of the moduli space $\mathcal{B}(-1,10)$ of a rank 2 stable bundles with the given Chern classes.
We also show that Hartshorne's conditions on a sequence $\mathcal{X}$ of 10 integers
are sufficient and necessary for the existence of a stable rank 2 bundle with odd determinant and spectrum $\mathcal{X}$. Furthermore, we prove that the sequence of 
integers $\{-2^{n-1},-1,0,1^{n-1}\}$ for $ n\geq4$ is realized as the spectrum of a stable rank 2 bundle $\EE$ of odd determinant by computing the minimal generators of its Rao module.
\end{abstract}

\tableofcontents

\section{Introduction}

There are three forms to classifying rank 2 vector bundles on the projective space $\mathbb{P}^3$. The first form is due to Horrocks since he proofs in \cite{Ho64} that every vector bundle on $\mathbb{P}^3$ is given as the cohomology of a \textit{monad} whose terms are summands of line bundles. Another way to study stable rank 2 vector bundles on the projective space is in terms of their \textit{spectra}, as originally introduced by Barth and Elencwajg in \cite{BE78}
Finally, we can also classify rank 2 vector bundles on $\mathbb{P}^3$ up to flat deformation by analyzing the irreducible components of the moduli space of stable vector bundles with Chern classes fixed. It is important to comment that there is no straightforward relation between these three classification schemes: a monad can be associated with two vector bundles with different spectra; a single spectrum can be associated with different monads; and an irreducible component of the moduli space may contain points corresponding to different spectra and monads.

Let $\mathcal{B}(e,m)$ denote the moduli space of stable rank 2 bundles on $\mathbb{P}^3$ with Chern classes $c_1=e$ and $c_2=m$. Up to normalization, it is enough to consider $e=-1,0$, and when $e=-1$ the integer $m$ must be even.

There is a complete classification of $\mathcal{B}(0,m)$ for $c_2\leq5$, more precisely: in \cite[Section 5.3]{HR91} it is shown that $\mathcal{B}(0,1)$ and $\mathcal{B}(0,2)$ are irreducible while in \cite{ES} and \cite{C} is proved 
that $\mathcal{B}(0,3)$ and $\mathcal{B}(0,4)$ have two irreducible components, respectively. The complete description of $\mathcal{B}(0,5)$ was obtained in \cite{AJTT} by proving that this moduli space has exactly three irreducible components.

When $e=-1$ we found fewer results in the literature: Hartshorne and Sols \cite{HS81} and Manolache \cite{M81} proved that $\mathcal{B}(-1,2)$ is irreducible; Banica and Manolache showed in \cite{BM85} that $\mathcal{B}(-1,4)$ has two irreducible components. More recently, the authors proved in \cite{MF2021} that both the moduli spaces $\mathcal{B}(-1,6)$ and $\mathcal{B}(-1,8)$ have at least four irreducible components.

The goal of the present paper is to list all \textit{positive} minimal Horrocks monads (see Section 2.2 for precise definitions) whose cohomology is a stable rank 2 bundle with odd determinant and $c_2=10$, to determine which ones exist, and to study the irreducible components of $\mathcal{B}(-1,10)$. We prove that, except for one case (namely $\mathcal{X}_8^{10}$, see Section 2.3), all possible spectra of Table \ref{spectra} are realized as the spectrum of a stable rank 2 bundle given by cohomology of a positive minimal monad. The remaining case is realized as the spectrum of a stable rank 2 bundle given as cohomology of a \textit{negative} minimal monad and we present it explicitly. In fact, we thank Nicolae Manolache for warning us about a mistake in our previous paper \cite{MF2021} regarding the spectrum $\mathcal{X}_6^{8}$; the same argument he indicated allowed us to correct our mistake and study the spectrum $\mathcal{X}_8^{10}$.

The paper is organized as follows. In Section \ref{section2} we recall three forms of classifying rank $2$ stable bundles: the Hartshorne--Serre correspondence, minimal Horrocks monads, and their cohomology, and the spectrum of a rank 2 bundle. By applying the Hartshorne--Serre correspondence to certain specific cases we get the families of \textit{Hartshorne bundles} and \textit{Ein bundles}. We also recall the definitions of positive, non-negative, and homotopy-free minimal monads and we revisit \cite[Lemma 3]{MF2021} which characterizes the minimal Horrocks monads whose cohomology is a stable rank 2 bundles on $\mathbb{P}^3$ with odd determinant. We finish this section by listing all possible spectra for rank 2 bundles with odd determinant and $c_2=10$.

In Section \ref{section3}, we enunciate Theorem \ref{possibles} which is the main tool to produce positive minimal monads with a stable rank 2 bundles with odd determinant as cohomology, and we apply this result to list all possible such positive monads when $c_2=10$ and fixed spectrum. We also 
eliminate some of the possibilities of positive minimal Horrocks monads.  

The following section is devoted to providing a complete classification of the positive minimal Horrocks monads whose cohomology is a rank $2$ stable vector bundle with Chern classes $c_1=-1$ and $c_2=10$; they are summarized in Table \ref{table:c2=10}. To prove the existence of these positive minimal monads, it was necessary to apply three methods: Lemma \ref{lema48}, an explicit construction of some positive minimal monads using Macaulay2, and the Hartshorne--Serre correspondence. Section \ref{negative-monads} is devoted to showing that the sequence of integers $\mathcal{X}_8^{10}:=\{-2^4,-1,0,1^4\}$ is realized as the spectrum of a negative minimal Horrocks monad and thus all possible spectra in Table \ref{spectra} are realized as the spectrum of a stable rank 2 bundle on $\mathbb{P}^3$. We also explicitly present the negative minimal monad whose cohomology is a stable bundle with spectrum $\{-2^{n-1},-1,0,1^{n-1}\}$ for $ n\geq4$.

In Section \ref{section5} we list all homotopy-free, positive minimal monads and compute the dimensions of the induced families of stable rank 2 bundles. Using the inferior semi-continuity of the dimension of cohomology groups, we identify a new irreducible component of $\mathcal{B}(-1,10)$ and show that $\mathcal{B}(-1,10)$ has at least four irreducible components, see Theorem \ref{dimension}. 

\subsection*{Acknowledgments}
MJ is supported by the CNPQ grant number 305601/2022-9, the FAPESP Thematic Project 2018/21391-1, and the FAPESP-ANR project 2021/04065-6.

\section{Constructing stable rank 2 bundles on $\mathbb{P}^3$}\label{section2}
\subsection{Hartshorne--Serre correspondence}\label{sub1}
The Hartshorne--Serre correspondence provides a tool to construct stable rank 2 bundles on 
$\mathbb{P}^3$ which establishes a relation between locally complete intersection subscheme 
$X\subset\mathbb{P}^3$ of codimension 2 and rank 2 bundles on $\mathbb{P}^3$ which has a section whose zero set is $X$. 
 More generally, 
\begin{Theorem}{\cite[Theorem 4.1]{H80}}\label{HS}
Fix an integer $c_1$. There is a one-to-one correspondence between:
\begin{enumerate}
    \item pairs $(\mathcal{E},s)$ where $\mathcal{E}$ is a rank $2$ reflexive sheaf on $\mathbb{P}^3$ with 
    $c_1(\mathcal{E})=c_1$ and $s\in H^0(\mathcal{E})$ is a global section whose zero-set has codimension $2$, and
    \item pairs $(X,\eta)$, where $X$ is a Cohen--Macaulay curve on $\mathbb{P}^3$, generically locally 
    complete intersection and $\eta\in H^0(\omega_X(4-c_1))$ is a global section which generates 
    the sheaf $\omega_X(4-c_1)$ except at finitely many points.
\end{enumerate}
\end{Theorem}
For this paper, it is sufficient to consider $\EE$ a rank 2 bundle on 
$\mathbb{P}^3$ and $X$ a local complete intersection of codimension 2. The section $s$ in Theorem \ref{HS} 
induces the exact sequence 
\begin{equation}\label{eq-HS}
0\longrightarrow\op3\longrightarrow\mathcal{E}\longrightarrow\mathcal{I}_X(k)\longrightarrow0,
\end{equation}
where $c_1(\mathcal{E})=k$ and $c_2(\mathcal{E})=\deg X$.

As a naive way to apply the correspondence in Theorem \ref{HS}, we remember the two families of 
vector bundles on $\mathbb{P}^3$ with $c_1=-1$ found in the literature: the first family of 
bundles $\{\EE_t\}_{t\geq2}$ \cite[Example 4.3.2]{H78} where each bundle $\EE_t$ of this family is 
obtained  applying the Hartshorne--Serre correspondence being $c_1=3$, the curve $X$ a disjoint union 
of $r$ conics $X_i$ on $\mathbb{P}^3$ and $\eta$ a section of $H^0(\omega_X(1))$. We get the exact sequence
\begin{equation*}
0\longrightarrow\op3\longrightarrow\mathcal{E}_t(2)\longrightarrow\mathcal{I}_X(3)\longrightarrow0,
\end{equation*}
with $c_1(\EE_t)=-1$ and $c_2(\EE_t)=2r-2$. The family of vector bundles so obtained is called 
\textit{Hartshorne family} and a bundle of this family is called \textit{Hartshorne bundle}. A vector 
bundle $E$ of the second family is obtained by 
considering a curve $X=X_1\cup X_2$ where $X_j$ is a complete intersection curve of bidegree $(n-b_i,n+b_i+1)$
for integers $b_j, i=1,2$. We have $\omega_X\simeq\mathcal{O}_X(2n+1)$ and 
from Theorem \ref{HS} we get the exact sequence
$$\op3\longrightarrow E'\longrightarrow\mathcal{I}_X(2n+1)\longrightarrow0$$ 
where $E=E'(-n-1)$ is such that $c_1(E)=-1$ and $c_2(E)=n^2-b_1^2-b_2^2+n-b_1-b_2$. We will call 
a vector bundle of this family of \textit{Ein bundle}.

\subsection{Minimal Horrocks Monads}

Another way to build bundles of rank $2$ on $\mathbb{P}^3$ is due to Horrocks \cite{Ho64}. Recall that a 
\textit{monad} on $\mathbb{P}^3$ is a complex

\begin{equation}\label{eq2}
\mathbf{M}:\ \ \  \ \mathcal{C}\stackrel{\alpha}{\longrightarrow}\mathcal{B}\stackrel{\beta}{\longrightarrow}\mathcal{A},
\end{equation}
of vector bundles on $\mathbb{P}^3$ such that the map $\alpha$ is injective and $\beta$ is surjective. The sheaf 
$\mathcal{E}:=\ker \beta/\Image\alpha$ is the \textit{cohomology of the monad} $\mathbf{M}$ of \eqref{eq2}. Let us assume that 
the morphism $\beta$ in \eqref{eq2} is locally left invertible, and so the cohomology of the 
monad $\EE$ of $\mathbf{M}$ is a vector bundle. Horrocks in 
\cite{Ho64} proves that all vector bundles on $\mathbb{P}^n$ can be obtained as cohomology of a 
monad of the form \eqref{eq2}, where the vector bundles $\mathcal{C}, \mathcal{B}$ 
and $\mathcal{A}$ are sums of line bundles.

In this paper, we 
always assume that the morphism $\alpha$ is locally left invertible, so that $\mathcal{E}$ is a vector 
bundle. A monad $\mathbf{M}$ of the form $\eqref{eq2}$ is called \textit{minimal} if the entries of the 
associated matrices to the maps $\alpha$ and $\beta$, as homogeneous forms, are not nonzero 
constants or, in other words, if no direct summand of $\mathcal{A}$ is the 
image of a line subbundle of $\mathcal{B}$ and if no direct summand of $\mathcal{C}$ goes 
onto a direct summand of $\mathcal{B}$. 
In addition, $\mathbf{M}$ is said to be
\textit{homotopy free} if
$$ \Hom(\mathcal{B},\mathcal{C})=\Hom(\mathcal{A},\mathcal{B})=0. $$
By considering $\EE$ be a stable rank 2 bundle on $\mathbb{P}^3$, there is an unique isomorphism 
$f:\EE\longrightarrow\EE^{*}(-1)$ with twisted sympletic structure that is $f^{*}(-1)=-f$. In \cite[Lemma 3]{MF2021}, it was shown 
the following characterization of a minimal monad whose cohomology is a stable rank 2 
bundle on $\mathbb{P}^3$ with $c_1=-1$, compare with \cite[Theorem 2.3]{MR2010} and \cite[Proposition 2.2]{R84}.


\begin{Lemma}\label{lema3}
Given a stable rank 2 bundle $\mathcal{E}$ on $\mathbb{P}^3$ with $c_1(\mathcal{E})=-1$, we can find two 
sequences of integers $\boldsymbol{a}=(a_1,\dots,a_s)$ and $\boldsymbol{b}=(b_1,\dots,b_{s+1})$ such 
that $\mathcal{E}$ is the cohomology of a monad of the form
\begin{equation}\label{eq6}
\bigoplus_{i=1}^{s}\op3(-a_i-1) \stackrel{\alpha}{\longrightarrow} \bigoplus_{j=1}^{s+1}
\Big(\op3(b_j)\oplus\op3(-b_j-1)\Big) \stackrel{\beta}{\longrightarrow} \bigoplus_{i=1}^{s}\op3(a_i),
\end{equation}
where we order the tuples $\boldsymbol{a}$ and $\boldsymbol{b}$ as $a_1\le \cdots \le a_s$ and 
$0\le b_1\le \cdots \le b_{s+1}$. 
In addition, 
\begin{equation}\label{eq1}
c_2(\mathcal{E})=\displaystyle\sum_{i=1}^{s}a_i(a_i+1)-\displaystyle\sum_{j=1}^{s+1}b_j(b_j+1).
\end{equation}
\end{Lemma}

Let us consider $\mathcal{A}=\bigoplus_{i=1}^{s}\op3(a_i)$ and $\mathcal{B}=\bigoplus_{j=1}^{s+1}
\Big(\op3(b_j)\oplus\op3(-b_j-1)\Big)$. Following the notation of Lemma \ref{lema3}, let 
$M:=H^1_*(\mathcal{E})$ and we suppose that $M$ has a minimal free presentation of the form
$$ \cdots\rightarrow F_1\rightarrow F_0\rightarrow M\rightarrow0.$$
It is known  that (see \cite[Proposition 2.2]{R84}) $\rank F_0=s$ and $\rank F_1=2s+2$ and conversely if $M$ 
has a free minimal resolution as above then with an argument similar to Hartshorne and Rao in 
\cite[Proposition 3.2]{HR91} it is shown that $\mathcal{E}$ is given as cohomology of a minimal 
monad of the form \eqref{eq6} with $F_0=H_*^0(\mathcal{A})$ and $F_1=H_*^0(\mathcal{B})$ and thus 
$\mathcal{A}, \mathcal{B}$ are the sheafifications of $F_0$ and $F_1$, respectively.\\
The minimal monad as in \eqref{eq6} is also called \textit{minimal Horrocks monad}. If we take the sequence 
$\boldsymbol{a}=(1,\ldots,1)$, then $b_j=0$ for $j=1,\ldots,s+1$ and this monad 
is called \textit{Hartshorne monad} whose cohomology is a vector bundle of the first family 
commented in Subsection \ref{sub1}. Another special monad in the literature is 
obtained by taking $s=1$ with $a_1>b_2\geq b_1\geq0$, 
this one is called \textit{Ein monad} whose cohomology is a stable bundle $\mathcal{F}$ with $c_1=-1$ and 
$c_2=a_1^2-b_1^2-b_2^2+a_1-b_1-b_2$,  this stable bundle belongs to the second family cited in 
Subsection \ref{sub1}.

Recall that a monad as in display \eqref{eq6} is called \textit{positive} if all summands of its right-hand term 
have a positive degree; that is $a_i>0$ for each $i=1,\dots,s$. If all summands of its right-hand term have 
a non-negative degree, then the monad is said to be \textit{non-negative}. If any of the summands in the 
right-hand term is negative, then the monad is said to be \textit{negative}.
\subsection{Spectrum of a vector bundle}

The \textit{spectrum} of a vector bundle $\mathcal{E}$ on $\mathbb{P}^3$ is a sequence of integers that 
encodes partial information about the cohomology modules $H^1_*(\mathcal{E})$ and 
$H^2_*(\mathcal{E})$. It was defined by Barth and Elencwajg \cite{BE78}, and 
Hartshorne extended this notion to rank 2 reflexive sheaves, see \cite[Theorem 7.1, p.151]{H80}. Following 
\cite[Theorem 7.1, pg. 151]{H80}, let $\mathcal{E}$ be a stable rank $2$
vector bundle on $\mathbb{P}^3$ with 
$c_1(\mathcal{E})=-1$; set $n:=c_2(\mathcal{E})$, and recall that this is necessarily 
even. The \textit{spectrum} of $\mathcal{E}$ is an unique sequence of 
integers $\X(\mathcal{E})=\{k_1,k_2,\cdots,k_n\}$ satisfying the 
following properties
\begin{enumerate}[label=(\alph*)]
\item \label{itemb1} $h^1(\mathbb{P}^3,\mathcal{E}(l)) = h^0(\mathbb{P}^1,\mathcal{H}(l+1))$ for $l\leq-1$ and
\item $h^2(\mathbb{P}^3,\mathcal{E}(l)) = h^1(\mathbb{P}^1,\mathcal{H}(l+1))$ for $l\geq-2$,
\end{enumerate}
where $\mathcal{H}=\bigoplus\op1(k_i)$. The spectrum $\mathcal{X}$ of a stable 
bundle $\EE$ satisfies the following properties which we 
organize in the form of a lemma and whose proof can be found in \cite[Section 7]{H80}
\begin{Lemma}\label{properties}
Let $\EE$ be a stable rank 2 bundle on $\mathbb{P}^3$ with $c_1(\EE)=-1, c_2(\EE)=n$. Then
\begin{enumerate}[label=\textbf{C.\arabic*}]
\item \label{itema1}(Symmetry)$\{k_i\}=\{-k_i-1\}$;
\item\label{itema2} (Connectness) Any integer $k$ between two integers of $\X$ also belongs to spectrum $\X$;
\item\label{itema3} If $k=\max\{-k_i\}$ and there is an integer $u$ with $-k\leq u\leq-2$ that 
occurs just once in $\X$, then each $k_i\in\X$ with $-k\leq k_i\leq u$ occurs exactly once, see \cite[Proposition 5.1]{H82}.
\end{enumerate}

\end{Lemma}
Before we list all the possibilities of spectra with $n=10$ integers we remember the notion 
of order between two spectra of same length: if $\mathcal{X}^n=\{k_1,k_2,\cdots,k_{n}\}$ and 
${\mathcal{S}}^n=\{k'_1,k'_2,\cdots,k'_{n}\}$ are spectra with $n$ integers in non-descending order, then
we say $\mathcal{X}^n>{\mathcal{S}}^n$ provided the left-most nonzero integer $k_i-k'_i$ 
is positive. Let us denote $r_j:=\{-j-1,j\}$ and $r_jr_i:=\{-j-1, -i-1, i, j\}$.

By applying Lemma \ref{properties}, we list on the Table \ref{spectra} all spectrum 
possibilities of a stable rank 2 vector bundle $\EE$ with $c_1(\EE)=-1$ and $c_2(\EE)=10$. On 
Table \ref{spectra} the first column stands the power of $r_0$ on the spectrum while 
the spectra in the lines of one are listed in decreasing order, from left to right.

\begin{table}[ht]
\begin{tabular}{|c|p{10.0cm}|}\hline
 $i$ &Spectrum\\
\hline
5& $\X_1^{10}=\{r_0^5\}$\\
\hline
4& $\X_2^{10}=\{r_0^4r_1\}$\\
\hline
3 & $\X_3^{10}=\{r_0^3r_1^2\}, \X_4^{10}=\{r_0^3r_1r_2\}$\\
\hline
3& $\X_5^{10}=\{r_0^2r_1^3\}, \X_6^{10}=\{r_0^2r_1^2r_2\}, \X_7^{10}=\{r_0^2r_1r_2r_3\}$\\
\hline
1&$\X_8^{10}=\{r_0r_1^4\}, \X_9^{10}=\{r_0r_1^3r_2\}, \X_{10}^{10}=\{r_0r_1^2r_2^2\}, \X_{11}^{10}=\{r_0r_1^2r_2r_3\}, 
\X_{12}^{10}=\{r_0r_1r_2r_3r_4\}$\\
\hline
\end{tabular}
\medskip
\caption{Possible spectra for stable rank $2$ vector bundles with $c_2(\mathcal{E})=10$.}
\label{spectra}
\end{table}



\section{Possible positive minimal Horrocks monads}\label{section3}

Let $\EE$ be a stable rank $2$ vector bundle on $\mathbb{P}^3$ with $c_1=-1$ and $c_2$ a positive 
integer. From the properties of spectrum follows that the spectrum $\X^{c_2}$ of $\EE$ can be written as
\begin{equation}\label{eq3}
\mathcal{X}^{c_2}=\{{(-k-1)}^{s(k)}...,0^{s(0)},\cdots, k^{s(k)}\}.
\end{equation}
If we denote $A=\mathbf{k}[X_0,X_1,X_2, X_3]$ to be the ring of homogeneous polynomials, then 
$M:=H_*^1(\EE)=\displaystyle\bigoplus_{l\in\mathbb{Z}}H^1(\EE(l))$ is called \textit{first cohomology module} 
of $\EE$ which is also a module over $A$. For each integer $l$ let us consider
$$\begin{array}{c}
m_l=h^1(\EE(l))=\dim M_l\\
\rho(l)=\dim\left[H^1(\EE(l))/\op3(1)\otimes H^1(\EE(l-1))\right],
\end{array}$$
that is, $\rho(l)$ is the number of minimal generators for $M$ in degree $l$ .

Fixed a possible spectrum in Table \ref{spectra}, we remember \cite[Theorem 8]{MF2021} which is applied to list all possible positive minimal Horrocks monads having $\EE$ as cohomology. 
\begin{Theorem}\label{possibles}
Let $\mathcal{E}$ be a stable rank 2 bundle on $\mathbb{P}^3$ with $c_1(\mathcal{E})=-1$ 
and spectrum as in display \eqref{eq3}. 
We have, 
\begin{equation}
\rho(-k-1)=m_{-k-1}=s(k)
\end{equation}
and for $0\leq i < k$,
\begin{equation}\label{ineq1}
s(i)-2\displaystyle\sum_{j\geq i+1}s(j)\leq\rho(-i-1)\leq s(i)-1.
\end{equation}
\end{Theorem}
From Table \ref{spectra} and Theorem \ref{possibles} we easily determine the existence of stable 
rank 2 bundles with spectrum $\X_1^{10}$ and $\X_{12}^{10}$, respectively. For the spectrum $\X_1^{10}$ 
we have $k=0$ and $\rho(-1)=5$, hence we get the Hartshorne monad with $s=5$ and $\mathbf{b}=(0,0,0,0,0,0)$
$$\op3(-2)^{\oplus5}\longrightarrow\op3^{\oplus5}\oplus\op3(-1)^{\oplus5}\longrightarrow\op3(1)^{\oplus5},$$
which exists and its cohomology is a stable rank $2$ bundle $\EE$ with spectrum $\X_1^{10}$. Now if we look 
to the spectrum $\X_{12}^{10}=\{r_0r_1r_2r_3r_4\}$, then we obtain 
$k=4, \rho(-5)=1$ and $\rho(-i-1)=0$, for $i=0, 1, 2, 3$. According to the equation of $c_2$ with $s=1$ 
follows that $b_1^2+b_2^2+b_1+b_2=20$ having as solution $b_1=4$ and $b_2=0$. In this form we 
get the Ein monad ($s=1, \mathbf{b}=(4,0)$ c. f. \cite{Ein88}) that exists and we write
$$\op3(-6)\longrightarrow\op3(4)\oplus\op3\oplus\op3(-1)\oplus\op3(-5)\longrightarrow\op3(5),$$ 
whose its cohomology is a stable rank $2$ bundle on $\mathbb{P}^3$ with spectrum $\X_{12}^{10}$. 

For each of the other spectra of Table \ref{spectra}, there is more than one possibility of minimal 
Horrocks monad and thus we will list all these possibilities in Table \ref{terms}.
\begin{table}[ht]
    \centering
    \begin{tabular}{|c|c|c|c|c|}
    \hline
     Spectrum & $k$& $\rho(-k-1)$& $i$& $\rho(-i-1)\in$\\
     \hline
     $\X_2^{10}$& 1& 1& 0& $\{2, 3\}$\\
     \hline
   $\X_3^{10}$& 1& 2& 0& $\{0, 1, 2\}$\\
   \hline
   \multirow{2}{*}{$\X_4^{10}$}&\multirow{2}{*}{2}&\multirow{2}{*}{1}& 0 &$\{0, 1, 2\}$\\
   \cline{4-5}
   &&& 1& $\{0\}$\\
   \hline
   $\X_5^{10}$& 1& 3& 0& $\{0, 1\}$\\
   \hline
   \multirow{2}{*}{$\X_6^{10}$}&\multirow{2}{*}{2}&\multirow{2}{*}{1}& 0 &$\{0, 1\}$\\
   \cline{4-5}
   &&& 1& $\{0, 1\}$\\
   \hline
   \multirow{2}{*}{$\X_7^{10}$}&\multirow{2}{*}{3}&\multirow{2}{*}{1}& 0&$\{0, 1\}$\\
   \cline{4-5}
   &&&$\{1, 2\}$&$\{0\}$\\
   \hline
   $\X_8^{10}$& 1& 4& 0& $\{0\}$\\
   \hline
   \multirow{2}{*}{$\X_9^{10}$}&\multirow{2}{*}{2}&\multirow{2}{*}{1}& 0 &$\{0\}$\\
   \cline{4-5}
   &&& 1& $\{1,2\}$\\
   \hline
   \multirow{2}{*}{$\X_{10}^{10}$}&\multirow{2}{*}{2}&\multirow{2}{*}{2}& 0 &$\{0\}$\\
   \cline{4-5}
   &&& 1& $\{0, 1\}$\\
   \hline
      \multirow{3}{*}{$\X_{11}^{10}$}&\multirow{3}{*}{3}&\multirow{3}{*}{1}& 0&$\{0\}$\\
   \cline{4-5}
&&&1&$\{0, 1\}$\\
\cline{4-5}
&&&2&$\{0\}$\\
   \hline
    \end{tabular}
    \caption{Number of minimal generators of $H_*^1(\EE)$}
    \label{terms}
\end{table}

Our goal is to list all possible positive minimal Horrocks monads with a fixed spectrum. For this, we 
take each possibility of number of minimal generators of $H_*^1(\EE)$ listed in Table 
\ref{terms} (that is for each sequence of integers $\mathbf{a}$) and we substitute on the 
equation \eqref{eq1} to find its solutions which are the sequences of integers 
$\mathbf{b}$. With the desire to simplify the notation, we denote
$$B_s:=\displaystyle\sum_{j=1}^{s+1}b_j(b_j+1).$$ 
We list all the possibilities of positive minimal Horrocks monads for each fixed spectrum c.f. Table \ref{terms}.

\vspace{1cm}
$\left(\mathbf{\X_2^{10}}, \rho(-2)=1\right)$
\begin{itemize}
    \item$\rho(-1)=2$ imply $\boldsymbol{a}=(2,1,1)$ and
    $$B_3=0\Leftrightarrow \boldsymbol{b}=(0,0,0,0).$$
    \item $\rho(-1)=3$ imply $\boldsymbol{a}=(2,1,1,1)$ and
    $$B_4=2\Leftrightarrow \boldsymbol{b}=(1,0,0,0,0).$$
\end{itemize}
$\left(\mathbf{\X_3^{10}}, \rho(-2)=2\right)$

\begin{itemize}
    \item$\rho(-1)=0$ imply $\boldsymbol{a}=(2,2)$ and
    $$B_2=2\Leftrightarrow \boldsymbol{b}=(1,0,0).$$
    \item $\rho(-1)=1$ imply $\boldsymbol{a}=(2,2,1)$ and
    $$B_3=4\Leftrightarrow \boldsymbol{b}=(1,1,0,0).$$
    \item $\rho(-1)=2$ imply $\boldsymbol{a}=(2,2,1,1)$ and
    $$B_4=6\Leftrightarrow \boldsymbol{b}=(1,1,1,0,0) \mbox{ or } \textcolor{red}{\boldsymbol{b}=(2,0,0,0,0)}.$$
\end{itemize}
$\left(\mathbf{\X_4^{10}}, \rho(-3)=1\right)$

\begin{itemize}
    \item$\rho(-1)=0$ imply $\boldsymbol{a}=(3)$ and
    $$B_1=2\Leftrightarrow \boldsymbol{b}=(1,0).$$
    \item $\rho(-1)=1$ imply $\boldsymbol{a}=(3,1)$ and
    $$B_2=4\Leftrightarrow \boldsymbol{b}=(1,1,0).$$
    \item $\rho(-1)=2$ imply $\boldsymbol{a}=(3,1,1)$ and
    $$B_3=6\Leftrightarrow \boldsymbol{b}=(1,1,1,0) \mbox{ or } \boldsymbol{b}=(2,0,0,0).$$
\end{itemize}

$\left(\mathbf{\X_5^{10}}, \rho(-2)=3\right)$
\begin{itemize}
    \item $\rho(-1)=0$ imply $\boldsymbol{a}=(2,2,2)$ and
    $$B_3=8\Leftrightarrow \boldsymbol{b}=(1,1,1,1)\mbox{ or } \textcolor{red}{\boldsymbol{b}=(2,1,0,0)}.$$
    \item $\rho(-1)=1$ imply $\boldsymbol{a}=(2,2,2,1)$ and
    $$B_4=10\Leftrightarrow \boldsymbol{b}=(1,1,1,1,1) \mbox{ or } \textcolor{red}{\boldsymbol{b}=(2,1,1,0,0)}.$$
\end{itemize}
$\left(\mathbf{\X_6^{10}}, \rho(-3)=1\right)$
\begin{itemize}
    \item $\rho(-1)=\rho(-2)=0$ imply $\boldsymbol{a}=(3)$ and
    $$B_1=2\Leftrightarrow \boldsymbol{b}=(1,0).$$
    \item $\rho(-1)=1$ and $\rho(-2)=0$ imply $\boldsymbol{a}=(3,1)$ and
    $$B_2=4\Leftrightarrow \boldsymbol{b}=(1,1,0).$$
    \item $\rho(-1)=0$ and $\rho(-2)=1$ imply $\boldsymbol{a}=(3,2)$ and
    $$B_2=8\Leftrightarrow \boldsymbol{b}=(2,1,0).$$
    \item $\rho(-1)=\rho(-2)=1$ imply $\boldsymbol{a}=(3,2,1)$ and
    $$B_3=10\Leftrightarrow \boldsymbol{b}=(2,1,1,0).$$
    
\end{itemize}
$\left(\mathbf{\X_7^{10}}, \rho(-4)=1\right)$
\begin{itemize}
    \item $\rho(-1)=0$ imply $\boldsymbol{a}=(4)$ and $B_1=10$ which has no solution.
    \item $\rho(-1)=1$ imply $\boldsymbol{a}=(4,1)$ and
    $$B_2=12\Leftrightarrow \boldsymbol{b}=(2,2,0)\mbox{ or } \boldsymbol{b}=(3,0,0).$$
\end{itemize}
$\left(\mathbf{\X_8^{10}}, \rho(-2)=4\right)$
\begin{itemize}
    \item $\rho(-1)=0$ imply $\boldsymbol{a}=(2,2,2,2)$ and
    $$
B_4=14\Leftrightarrow \textcolor{red}{\boldsymbol{b}=(3,1,0,0,0)}\mbox{ or }
\textcolor{red}{\boldsymbol{b}=(2,2,1,0,0)} \mbox{ or } \textcolor{red}{\boldsymbol{b}=(2,1,1,1,1)}.
 $$
\end{itemize}
$\left(\mathbf{\X_9^{10}}, \rho(-3)=1\right)$
\begin{itemize}
  \item $\rho(-2)=1$ imply $\boldsymbol{a}=(3,2)$ and
  $$B_2=8\Leftrightarrow \boldsymbol{b}=(2, 1, 0).$$
  \item $\rho(-2)=2$ imply $\boldsymbol{a}=(3,2,2)$ and
  $$B_3=14\Leftrightarrow \boldsymbol{b}=(2,2,1,0)\mbox{ or } \textcolor{red}{\boldsymbol{b}=(3,1,0,0)}.$$
\end{itemize}
$\left(\mathbf{\X_{10}^{10}}, \rho(-3)=2\right)$
\begin{itemize}
    \item $\rho(-2)=0$ imply $\boldsymbol{a}=(3,3)$ and
    $$B_2=14\Leftrightarrow \boldsymbol{b}=(2,2, 1)\mbox{ or } \textcolor{red}{\boldsymbol{b}=(3,1,0)}.$$
   \item $\rho(-2)=1$ imply $\boldsymbol{a}=(3,3,2)$ and
    $$B_3=20\Leftrightarrow \boldsymbol{b}=(2,2,2,1)\mbox{ or } \textcolor{red}{\boldsymbol{b}=(3,2,1,0)}
    \mbox{ or } \textcolor{red}{\boldsymbol{b}=(4,0,0,0)}.$$
\end{itemize}
$\left(\mathbf{\X_{11}^{10}}, \rho(-4)=1\right)$
\begin{itemize}
    \item $\rho(-2)=0$ imply $\boldsymbol{a}=(4)$ and $B_1=10$ which has no solution.
    \item $\rho(-2)=1$ imply $\boldsymbol{a}=(4,2)$ and
    $$B_2=16\Leftrightarrow \boldsymbol{b}=(3,1,1).$$
\end{itemize}
\begin{prop}\label{eli}
If one of the 10 possibilities of positive minimal Horrocks monads in red listed exist then its cohomology is not a stable rank 2 bundle.
\end{prop}
\begin{proof}
Suppose, for example, that the positive minimal Horrocks 
monad $\mathbf{M}$ in which $\mathbf{a}=(2,2,1,1)$ 
and $\mathbf{b}=(2,0,0,0,0)$ exists (the other possibilities can be treated analogously) and its cohomology is a bundle $\EE$. We can write it as
$$ \op3^{\oplus2}(-2){\oplus}\op3^{\oplus2}(-3) \stackrel{\alpha}{\longrightarrow} \op3(2){\oplus}\op3^{\oplus4}{\oplus}\op3^{\oplus4}(-1){\oplus}\op3(-3) \stackrel{\beta}{\longrightarrow} \op3^{\oplus2}(2){\oplus}\op3^{\oplus2}(1),$$
where $\alpha$ and $\beta$ are maps injective and surjective, respectively. But from the minimally 
of the monad $\mathbf{M}$ the first column of the map $\beta$ is zero which implies 
$\beta\circ \imath=0$, where $\imath$ denotes the inclusion of $\op3$ into the first 
summand of $\op3(2)\oplus\op3^{\oplus4}\oplus\op3^{\oplus4}(-1)\oplus\op3(-3)$. This means that 
$\imath\in H^0(\ker\beta)=H^0(\EE)$. Therefore, $\E$ is not stable.
\end{proof}

\begin{Lemma}\label{lemma-a1}
Let $a, b$ be positive integers such that $b<a$. Consider the possible positive  minimal Horrocks monad
\begin{equation}\label{monad1}
\op3(-b-1)\oplus \op3^{\oplus3}(-a-1) \stackrel{\alpha}{\longrightarrow}
\op3^{\oplus5}(b)\oplus \op3^{\oplus5}(-b-1) \stackrel{\beta}{\longrightarrow}
\op3(b)\oplus \op3^{\oplus3}(a).
\end{equation}
\begin{enumerate}
    \item If $6b+1\geq 4a$, then the possible of monad in \eqref{monad1} does not exist.
    \item For $6b+1\geq3a$, if the monad at in \eqref{monad1} exists then its cohomology is a vector bundle unstable.
\end{enumerate}
\end{Lemma}

\begin{proof}
From the minimally of the monad follows that $\op3^{\oplus5}(b)$ maps on $\op3^{\oplus3}(a)$, let us say 
$\phi:\op3^{\oplus5}(b)\longrightarrow \op3^{\oplus3}(a)$. The possibilities of the rank of 
$\phi$ are $1, 2, 3$. If $\rank\phi=j$ with $j=1, 2$ then $h^0((\Image\phi)(1)))\leq j{a+4\choose3}$ which implies $h^0(\Ker\phi(1))\geq5{b+4\choose3}-j{a+4\choose3}$. Therefore, 
$$h^0(\EE(1))\geq h^0(\Ker\phi(1))\geq5{b+4\choose3}-j{a+4\choose3}\geq10.$$
On the other hand, $h^0(\mathcal{I}_Y(1))=h^0(\EE(1))-1\geq9$ and this means that the curve $Y$ is contained in 
at least $9$ distinct planes, a contradiction.

Now we assume that $\rank\phi=3$. We can write
$$\beta=\left(
\begin{array}{cc}
    \beta' & 0 \\
   \beta''  & \phi
\end{array}
\right),$$
where $\beta':\op3(-b-1)^{\oplus5}\longrightarrow \op3(b)$ and $\beta'':\op3(-b-1)^{\oplus5}\longrightarrow \op3(a)^{\oplus3}$. Since 
$$\alpha=\Omega^{-1}\circ\beta^*(-1)=
\left(\begin{array}{cc}
    0 & -\phi^*(-1) \\
   \beta'^{*}(-1)  & \beta''^{*}(-1)
\end{array}
\right),
$$
follows that
$$\beta\alpha=\left(\begin{array}{cc}
    0 & -\beta'\phi^*(-1) \\
   \phi\beta'^{*}(-1)  & -\beta''\phi^*(-1)+\phi\beta''^{*}(-1)
\end{array}
\right).
$$
In particular, $\phi\beta'^{*}(-1)=0$ which implies $\Image\beta'^{*}(-1)\subset\ker\phi$ and we get the diagram
$$\xymatrix{
 \op3(-b-1)\ar[r]^{\beta'^{*}(-1)}& \op3^{\oplus5}(b) \ar[d]^{\phi}\ar[r]\ar[d] & G=\coker\beta'^{*}(-1) \ar[d]^{\overline{\phi}} \\
 & \op3^{\oplus3}(a) \ar[r] & \op3^{\oplus3}(a).\\	
}$$
From $\beta'$ is surjective and $\rank\phi=3$ the sheaf $G$ is locally free of rank $4$ 
and $\ker\overline{\phi}\simeq\op3(k)$. We 
observe that $H^0(G(-b-1))=0$ and analyzing the diagram
$$\xymatrix{
 \op3(k) \ar@{^{(}->}[r]& G\ar[rr]^{\overline{\phi}}\ar@{->>}[rd] && \op3(a)^{\oplus3}  \\
 &  & \Image\overline{\phi}\ar@{^{(}->}[ru]&\\	
}$$
we get $H^0(G(-k))\neq0$, thus $k\leq b$.  Also by the above diagram 
$$\mu(\mathcal{O}(k))\leq b<\frac{6b+1}{4}=\mu(G)<\mu(\Image\overline{\phi})\leq a.$$
This means that if the possibility of a positive minimal monad in \eqref{monad1} exists then 
$6b+1<4a$. Furthermore, we have $c_1(\Image\overline{\phi})=6b+1-k$ and so 
$k\geq6b+1-3a\geq0$ whenever $6b+1\geq3a$ which implies $H^0(E(-k))\neq0$ and $\EE$ is not stable.
\end{proof}
\begin{cor}\label{propa3}
If bundle $\EE$ has spectrum $\X_5^{10}$ and it is the cohomology of the positive minimal Horrocks monad 
$$ \op3(-2)\oplus\op3(-3)^{\oplus3} \stackrel{\alpha}{\longrightarrow}
\op3(1)^{\oplus5}\oplus\op3(-2)^{\oplus5} \stackrel{\beta}{\longrightarrow}
\op3(1)\oplus\op3(2)^{\oplus3}
$$
then $\EE$ cannot be stable.
\end{cor}
\begin{proof}
Take $b=1, a=2$ on Lemma \ref{lemma-a1}.
\end{proof}
\begin{prop}\label{nostable}
Let $a_2, a_1, b$ integers such that $a_2\geq b>a_1\geq0$ and $2b-a_2\geq0$. The vector bundle 
$\EE$ that is the cohomology of the monad 
$$\op3(-a_2-1)\oplus \op3(-a_1-1)^{\oplus g} \stackrel{\alpha}{\longrightarrow}
\op3(b)^{\oplus2}\oplus B \stackrel{\beta}{\longrightarrow}
\op3(a_2)\oplus \op3(a_1)^{\oplus g},
$$
is not stable, where $B$ is a $(2g+2)$-bundle .
\end{prop}
\begin{proof}
The display of the monad implies $H^0(\ker\beta(-p))\simeq H^0(\EE(-p))$ whenever $p\geq0$. By the minimality of the monad we have
$$\beta=\left(
\begin{array}{cc}
    \psi & \phi \\
   0  & \eta
\end{array}
\right),$$
where $\psi:\op3(b)^{\oplus2}\longrightarrow\op3(a_2)$. Since $\ker\psi$ is a reflexive sheaf of rank $1$ 
follows that $\ker\psi\simeq\op3(k)$ and we get the commutative diagram
$$\xymatrix{
 \op3(k)\ar@{^{(}->}[d]^{s}\ar[r]& \op3(b)^{\oplus2} \ar@{^{(}->}[d]\ar[r]^{\psi} & \op3(a_2)\ar@{^{(}->}[d] \\
 \ker\beta\ar[r]& \op3(b)^{\oplus2}\oplus B \ar[r]^{\beta} & \op3(a_2)\oplus \op3(a_1)^{\oplus g}.\\
}$$
Thus $s\in H^0(\ker\beta(-k))$ with $k=2b-c_1(\Image\psi)\geq2b-a_2\geq0$. 
Therefore, $s\in H^0(\ker\EE(-k)), k\geq0$ and so $\EE$ is not stable. 
\end{proof}
\begin{cor}\label{nonstable1}
If the possible positive minimal Horrocks monad 
$$ \op3(-5)\oplus\op3(-2) \stackrel{\alpha}{\longrightarrow}
\op3(2)^{\oplus2}\oplus\op3\oplus\op3(-1)\oplus\op3(-3)^{\oplus2} \stackrel{\beta}{\longrightarrow}
\op3(4)\oplus\op3(1)
$$
exists and $\EE$ is its cohomology with spectrum $\X_7^{10}$, then $\EE$ cannot be stable.
\end{cor}
\begin{proof}
Take $a_2=4, a_1=1, b=2$ and $g=1$ in Proposition \ref{nostable}.
\end{proof}
\begin{Lemma}\label{slope}
Let $E,F$ be torsion free semistable (stable) sheaves such that $\mu(E)>\mu(F)$ (resp. $\mu(E)\geq\mu(F)$). Then, $\Hom(E, F)=0$.
\end{Lemma}
\begin{proof}
Let $\varphi:E\rightarrow F$ be a morphism in $\Hom(E, F)$. If this one is injective then $\mu(E)\leq\mu(F)$ (resp. $\mu(E)<\mu(F)$) 
which is a contradiction. If $\varphi$ is not injective we consider the commutative diagram 
$$\xymatrix{
 \ker\varphi \ar[r]& E\ar[rr]^{\varphi}\ar[rd] && F  \\
 &  & \Image\varphi\ar@{^{(}->}[ru]&\\	
}$$
and we get $\mu(E)\leq\mu(\Image\varphi)\leq\mu(F)$ 
(resp. $\mu(E)<\mu(\Image\varphi)<\mu(F)$), which is also a contradiction. Thus $\Hom(E, F)=0$.
\end{proof}

\begin{prop}\label{nonstable2}
The cohomology of the monad
$$\op3(-4)\oplus \op3(-2)^{\oplus2} \stackrel{\alpha}{\longrightarrow}
\op3(1)^{\oplus3}\oplus B \stackrel{\beta}{\longrightarrow}
\op3(3)\oplus \op3(1)^{\oplus2}
$$
is not a stable bundle, where $B=\op3\oplus\op3(-1)\oplus\op3(-2)^{\oplus3}$.
\end{prop}
\begin{proof}
Let $\EE$ be the cohomology of this monad and $K=\ker\beta$. We suppose that the vector bundle $\EE$ is 
stable. By minimally of the monad we write
$$\beta=\left(
\begin{array}{cc}
    \psi & \eta \\
   0  & \phi
\end{array}
\right),$$
where $\psi:\op3(1)^{\oplus3}\rightarrow\op3(3)$. If we denote $P:=\ker\psi$ then we have the commutative diagram
$$\xymatrix{
 P\ar@{^{(}->}[d]\ar[r]& \op3(1)^{\oplus3} \ar@{^{(}->}[d]\ar[r]^{\psi} & \op3(3)\ar@{^{(}->}[d] \\
 K\ar[r]\ar[d]& \op3(1)^{\oplus3}\oplus B \ar[r]^{\beta} \ar@{^{(}->}[d]& \op3(3)\oplus \op3(1)^{\oplus2}\ar[d].\\
 K'\ar[r]&B\ar@{->>}[r]^{\phi}&\op3(1)^{\oplus2}.\\
}$$
Since $P\simeq\ker\{K\rightarrow K'\}$ and $K, K'$ are locally free sheaves 
follows that $P$ is a rank $2$ stable reflexive sheaf with $c_1 P=0$. On the other hand, Serre's duality implies
$$\begin{array}{c}
   \Hom(P,\op3(-4)\oplus\op3(-2)^{\oplus2}\simeq\ext^0(P,\op3(-4))\oplus2\cdot\ext^0(P,\op3(-2))   \\
   \simeq H^3(\mathbb{P}^3,P)^{*}\oplus 2\cdot H^3(\mathbb{P}^3,P(-2))^{*}\simeq H^0(P(-4))\oplus 2\cdot H^0(P(-2))=0.   
\end{array}$$
Since $0=c_1 P>c_1\EE=-1$, follows from Lemma \ref{slope} that
the map $P\rightarrow K\rightarrow\EE$ is zero and 
it can be extended to a map $\sigma:P\rightarrow \op3(-4)\oplus\op3(-2)^{\oplus2}$ 
$$\xymatrix{
 & P\ar[d] \ar@{-->}[ld]^{\sigma}\ar[rd]^{0}& \\
 \op3(-4)\oplus\op3(-2)^{\oplus2}\ar[r]& K \ar[r]& \EE\\
}$$
which is a contradiction. Therefore, $\EE$ is not stable.

\end{proof}
\section{Existence of positive minimal monads}\label{section4}
The goal of this section is to show the existence of positive minimal Horrocks monads whose cohomology 
is a stable rank 2 bundle on $\mathbb{P}^3$ with $c_2=10, c_1=-1$. With this wish, we 
remember the result below (see \cite[Lema 6]{MF2021}) which allows 
us to build up minimal monads for higher Chern classes. 

\begin{Lemma}\label{lema48}
Let $(E,\sigma)$ be pair consisting of a stable rank 2 vector bundle $E$ with $c_1(E)=-1$ and $c_2(E)=n$ 
and a section $\sigma\in H^0(E(r))$ with $r>0$ such that $X:=(\sigma)_0$ is curve. If $C$ is a 
complete intersection curve of type $(u,v)$ disjoint from $X$ satisfying $u+v=2r-1$, then there 
is a pair $(E',\sigma')$ consisting of a a stable rank 2 vector bundle $E$ with $c_1(E)=-1$ and 
$c_2(E)=n+uv$ and a section $\sigma'\in H^0(E'(r))$ such that $(\sigma')_0=X\sqcup C$. Moreover, 
if $E$ is the cohomology of a minimal monad of the form
$$ \mathbf{M}:\ \ \  \ \mathcal{C} \longrightarrow \mathcal{B} \longrightarrow \mathcal{A}, $$
then $E'$ is the cohomology of a minimal monad of the form
$$\mathbf{M}':\ \ \op3(-r)\oplus\mathcal{C} \longrightarrow \op3(r-1-u)\oplus\op3(r-1-v)
\oplus\mathcal{B} \longrightarrow \op3(r-1)\oplus\mathcal{A}.$$
\end{Lemma}

In order to apply Lemma \ref{lema48} we tabulate in Table \ref{c2<=8} all the positive minimal Horrocks monads whose 
cohomology is a stable rank 2 bundle on $\mathbb{P}^3$ with $c_2\leq8$ c.f. listed in \cite{MF2021}.  
	\begin{table}[ht]
		\begin{center}
	\begin{tabular}{ | l | c | c | c | } 
\hline
Spectrum & $\boldsymbol{b}$ &  $\boldsymbol{a}$ & Label\\ 
\hline
$\X_1^2=\{r_0\}$ & $(0,0)$ & $(1)$ & $\mathbf{M}_1$\\
\hline
$\X_1^4=\{r_0^2\}$ & $(0,0,0)$ & $(1,1)$ & $\mathbf{M}_2$\\
\hline
$\X_2^4=\{r_0r_1\}$ & $(1,0)$ & $(2)$ & $\mathbf{M}_3$\\
\hline
$\X_1^6=\{r_0^3\}$ & $(0,0,0,0)$ & $(1,1,1)$ & $\mathbf{M}_4$\\
\hline
$\X_2^6=\{r_0^2r_1\}$ & \makecell{$(0,0)$ \\ $(1,0,0)$} & \makecell{$(2)$ \\ $(2,1)$} & \makecell{$\mathbf{M}_5$ \\ $\mathbf{M}_6$}\\
\hline
$\X_3^6=\{r_0r_1^2\}$ &$(1,1,1)$ &$(2,2)$ &$\mathbf{M}_7$\\
\hline
$\X_4^6=\{r_0r_1r_2\}$ & $(2,0)$ & $(3)$ & $\mathbf{M}_8$\\
\hline
$\X^8_1=\{r_0^4\}$&$(0,0,0,0,0)$ &$(1,1,1,1)$ &  $\mathbf{M}_9$\\
\hline
$\X^8_2=\{r_0^3r_1\}$& \makecell{$(0,0,0)$\\ $(1,0,0,0)$}& \makecell{$(2,1)$\\ $(2,1,1)$} & \makecell{$\mathbf{M}_{10}$\\ $\mathbf{M}_{11}$}\\
\hline
$\X^8_3=\{r_0^2r_1^2\}$& \makecell{$(1,1,0)$\\ $(1,1,1,0)$}&  \makecell{$(2,2)$\\ $(2,2,1)$}& \makecell{$\mathbf{M}_{12}$\\ $\mathbf{M}_{13}$}\\
\hline
$\X^8_4=\{r_0^2r_1r_2\}$ &$(2,0,0)$& $(3,1)$&  $\mathbf{M}_{14}$ \\
\hline
$\X^8_5=\{r_0r_1^2r_2\}$&\makecell{$(1,1)$\\ $(2,1,1)$}& \makecell{$(3)$\\ $(3,2)$} & \makecell{$\mathbf{M}_{15}$\\ $\mathbf{M}_{16}$}\\
\hline
$\X^8_7=\{r_0r_1r_2r_3\}$ &$(3,0)$ &$(4)$ &$\mathbf{M}_{17}$\\
\hline
	\end{tabular}
\caption{Positive minimal Horrocks monads for $c_2\leq8$, see \cite{MF2021}.}
\label{c2<=8}
\end{center}
\end{table}

In Table \ref{table:c2=10} it will tabulate all positive minimal Horrocks monads 
whose cohomology is a stable bundle on $\mathbb{P}^3$ with $c_1=-1$ and $c_2=10$. Let us 
establish the following notation in Table: \\
\begin{enumerate}
    \item[(a)] The fourth column indicates the degree to which the global section exists for constructing the vector 
    bundle $\EE$, giving a curve $Y$ as a zero scheme.
    \item[(b)] The last column describe as we construct the curve $X$ by applying Lemma \ref{lema48}.
    \item [(c)] $P_n$ denotes a plane curve of degree $n$ while $C_{a,b}$ denotes a curve of bidegree $(a,b)$ on a smooth quadric.
\end{enumerate}

	\begin{table}[ht]
		\begin{center}
	\begin{tabular}{ | l | c | c | c |c|} 
\hline
Spectrum & $\boldsymbol{b}$ &  $\boldsymbol{a}$ &$r$ &Construction\\ 
\hline
$\X^{10}_1$ & $(0,0,0,0,0,0)$ & $(1,1,1,1,1)$ &$1$ & disjoint union of $6$ conics\\
\hline
\multirow{2}{4em}{$\X^{10}_2$}&$(0,0,0,0)$ &$(2,1,1)$& $3$&$\mathbf{M}_{2}, (3,2)$\\
&$(1,0,0,0,0)$& $(2,1,1,1)$&$3$&$\mathbf{M}_{4}, (4,1)$\\
\hline
\multirow{2}{4em}{$\X^{10}_3$}& $(1,0,0)$& $(2,2)$&$3$&$\mathbf{M}_{3}, (3,2)$\\
& $(1,1,0,0)$& $(2,2,1)$ &$3$&$\mathbf{M}_{6}, (4,1)$\\
& $(1,1,1,0,0)$& $(2,2,1,1)$ &$2$&$\mathbf{M}_{13}, (2,1)$\\
\hline
\multirow{4}{4em}{$\X^{10}_4$}& $(1,0)$& $(3)$&$2$&Ein\\
& $(1,1,0)$&$(3,1)$ &$2$&$\mathbf{M}_{15}, (2,1)$\\
& $\textcolor{blue}{(1,1,1,0)}$&$\textcolor{blue}{(3,1,1)}$ &-&-\\
& $(2,0,0,0)$&$(3,1,1)$ &$4$&$\mathbf{M}_{2}, (6,1)$\\
\hline
$\X^{10}_5$& $(1,1,1,1)$& $(2,2,2)$&$3$&$\mathbf{M}_{7}, (4,1)$\\
\hline
\multirow{4}{4em}{$\X^{10}_6$}& $(1,0)$& $(3)$&$2$&$C_{3,3}$\\
& $(1,1,0)$&$(3,1)$ & $2$ &$\mathbf{M}_{15}, (2,1)$\\
& $(2,1,0)$&$(3,2)$ &$3$&$\mathbf{M}_{8}, (4,1)$\\
& $(2,1,1,0)$&$(3,2,1)$ &$2$&$\mathbf{M}_{16}, (2,1)$\\
\hline
 \multirow{2}{4em}{$\X^{10}_7$}& $\textcolor{blue}{(2,2,0)}$& $\textcolor{blue}{(4,1)}$&-&-\\
 & $(3,0,0)$&$(4,1)$ & $2$ &$\mathbf{M}_{17}, (2,1)$\\
\hline
\multirow{2}{4em}{$\X^{10}_9$}& $(2,1,0)$& $(3,2)$&$3$&$\mathbf{M}_{8}, (4,1)$\\
 & $(2,2,1,0)$&$(3,2,2)$ &$1$& $P_2\cup P_3$ joined by one point\\
 \hline
 \multirow{2}{4em}{$\X^{10}_{10}$}& $(2,2,1)$& $(3,3)$&$1$& Explicit computation\\
 & $(2,2,2,1)$&$(3,3,2)$&$1$&$P_2\cup P_3$ joined by two points\\
 \hline
$\X^{10}_{11}$& $(3,1,1)$& $(4,2)$&$1$&Explicit computation\\
\hline
$\X_{10}^{12}$ & $(4,0)$ & $(5)$&$1$& $P_5$\\
\hline

	\end{tabular}
\caption{Positive minimal Horrocks monads for $c_2=10$. Was showed in Corollary \ref{nonstable1} and 
Proposition \ref{nonstable2} that if the possible monads in blue exist then its cohomologies are not stable bundles.}
\label{table:c2=10}
\end{center}
\end{table}
\newpage
\begin{prop}\label{prop3}
There are the possible positive minimal Horrocks monads with $\boldsymbol{b}=(3,1,1), \boldsymbol{a}=(4,2)$ 
and $\boldsymbol{b}=(2,2,1), \boldsymbol{a}=(3,3)$.
\end{prop}
\begin{proof}
To the prove the existence of the monad with $\boldsymbol{b}=(3,1,1)$ and  $\boldsymbol{a}=(4,2)$ it is enough to consider the maps 
$$\beta=\left(
\begin{array}{cccccc}
y& 0&w^3& z^6& 0& x^8\\
0& w& y& w^4&x^4&z^6\\ 
\end{array} \right)
\mbox{ and }\alpha=\left(
\begin{array}{cc}
-x^8&-z^6\\
-w^6&-x^4-yw^3\\
-z^6&0\\
w^3&y\\
0&w\\
y&0\\
\end{array}
\right),$$
where $\alpha$ is injective and $\beta$ and is surjective satisfying $\beta\circ\alpha=0$. Now if we have the other possible monad 
$\boldsymbol{b}=(2,2,1), \boldsymbol{a}=(3,3)$, then we take
$$\beta=\left(
\begin{array}{cccccc}
w& z&x^2& 0& x^6& y^6\\
z& x& w^2& w^5&y^6&0\\ 
\end{array} \right)
\mbox{ and }\alpha=\left(
\begin{array}{cc}
y^6-xzw^4&-x^2w^4\\
xw^5-x^3zw^2&y^6-x^4w^2\\
x^4z+xz^2w^2&x^5+x^2zw^2+w^5\\
-x^2&-w^2\\
-z&-x\\
-w&-z\\
\end{array}
\right).$$
It is easy to see that $\alpha$ is injective and $\beta$ is surjective such that $\beta\circ\alpha=0$. Therefore, 
the positive minimal Horrocks monad
$$\op3(-4)^{\oplus 2}\stackrel{\alpha}{\longrightarrow} \op3(2)^{\oplus 2}{\oplus}
\mathcal{O}(1){\oplus}\mathcal{O}(-2){\oplus}\op3(-3)^{\oplus 2} 
\stackrel{\beta}{\longrightarrow} \op3(3)^{\oplus 2}$$
indeed exists.
\end{proof}
Now, let $X$ be a locally complete intersection curve in $\mathbb{P}^3$ of degree $d$ which is a union of irreducible 
nonsingular curves meeting quasi-transversely. From \cite[Proposition 2.8]{HR91}, the sheaf 
$\mathcal{N}_X\otimes\omega_X(3)$ admits a general section which is nowhere vanishing and by 
considering $Y$ to be the multiplicity $2$ structure on $X$ defined by the associated exact sequence
\begin{equation}\label{exact11}
0\longrightarrow\mathcal{I}_Y\longrightarrow\mathcal{I}_X\longrightarrow\omega_X(3)\longrightarrow0,
\end{equation}
follows by Ferrand's Theorem \cite[Theorem 1.5]{H78} that $\omega_Y=\mathcal{O}_Y(-3)$. By Theorem
\ref{HS} there is a stable rank $2$ bundle $\EE(1)$ on $\mathbb{P}^3$ with $c_1(\EE)=-1, c_2(\EE)=2d$ such that the 
curve $(Y,\mathcal{O}_Y)$ is given as zero set of a section $s$ of $H^0(Y,\EE(1))$ and 
this section $s$ induces the exact sequence
\begin{equation}\label{exact12}
0\longrightarrow\op3(-1)\longrightarrow\EE\longrightarrow\mathcal{I}_Y\longrightarrow0.
\end{equation}
Since $c_1(\EE(1))=1$, if $X\neq\emptyset$, then the curve $Y$ of even degree does not 
lie on a hyperplane and so the vector bundle $\EE$ is stable. In this manner, we prove the following.
\begin{prop}\label{prop12}
Let $X$ be a locally complete intersection curve in $\mathbb{P}^3$ of degree $d$ which is a union of irreducible 
nonsingular curves meeting quasi-transversely where each connected component $C$ has $H^0(C,\mathcal{O}_C)=k$. Then 
there is a stable rank $2$ bundle $\mathcal{E}$ on $\mathbb{P}^3$ such 
that $c_1(\mathcal{E})=-1$ and $c_2(\mathcal{E})=2d$.  
\end{prop}
\begin{Theorem}\label{teo12}
The following $\mathbf{M}$ and $\mathbf{M'}$ possible positive minimal Horrocks monads  with 
$\boldsymbol{b}=(2,2,1,0), \boldsymbol{a}=(3,2,2)$ and $\boldsymbol{b}=(2,2,2,1), \boldsymbol{a}=(3,3,2)$, respectively, exist.
\end{Theorem}
\begin{proof}
Let $X$ be a locally complete intersection curve on $\mathbb{P}^3$ which is a 
union of irreducible nonsingular curves meeting quasi-transversely. From exact sequences 
\eqref{exact11} and \eqref{exact12} we get the following long exact sequence in cohomology
\begin{equation}\label{long-exact}
0\longrightarrow H^0(\omega_X(l+3))\longrightarrow H^1(\EE(l))\longrightarrow H^1(\mathcal{I}_X(l))\longrightarrow H^1(\omega_X(l+3))\longrightarrow\cdots
\end{equation}
for each integer $l\leq-1$ which implies
$$H^1(\EE(l))\simeq H^0(\omega_X(l+3)), l\leq-1.$$
If $X=X_1\cup X_2$ with $S=X_1\cap X_2$ containing $r$ collinear points, then we have (see \cite[Lemma 2.7]{HR91} the exact sequence
\begin{equation*}
0\rightarrow\omega_{X_1}\oplus\omega_{X_2}\rightarrow\omega_X\rightarrow\omega_S\rightarrow0
\end{equation*}
that implies 
\begin{equation}\label{dim}
h^0(\omega_X(m))=h^0(\omega_{X_1}(m))+h^0(\omega_{X_2}(m))+h^0(\omega_S(m))-\rank(\delta),
\end{equation}
where $\delta$ is the connecting homomorphism 
$$\delta:H^0(\omega_S(m))\longrightarrow H^1(\omega_{X_1}(m))\oplus H^1(\omega_{X_2}(m)),$$
such that 
$$\rank(\delta)=\left\{\begin{array}{cl}
  0   & \mbox{ if } m>0 \\
  1-m   & \mbox{ if } 2-r\leq m\leq0 \\
r   & \mbox{ if } m\leq1-r. \\
\end{array}
\right.
$$
If we consider the curve $X=P_2\cup P_3$ union of two plane curves $P_2$ and $P_3$ of degree $2$ and 
$3$, respectively, joined by one point $x_0$ then  $S=\{x_0\}$ and we have
$$\omega_{P_2}=\mathcal{O}_{P_2}(-1), \omega_{P_3}=\mathcal{O}_{P_3} \mbox{ and } h^0(\omega_{S}(m))=1, m\in\mathbb{Z}.$$
By Proposition \ref{prop12} we get a stable rank $2$ bundle $\EE$ on $\mathbb{P}^3$ such that 
$c_1(\EE)=-1$ and $c_2(\EE)=10$. From exact sequence \eqref{long-exact} and Equation \eqref{dim}
$$\begin{array}{l}
 h^1(\EE(-1))=h^0(\omega_X(2))= h^0(\mathcal{O}_{P_2}(1))+h^0(\mathcal{O}_{P_3}(2))+h^0(\omega_S(2))-0=10\\
 h^1(\EE(-2))=h^0(\omega_X(1))=h^0(\mathcal{O}_{P_2})+h^0(\mathcal{O}_{P_3}(1))+h^0(\omega_S(1))-0=5\\
 h^1(\EE(-3))=h^0(\omega_X)=h^0(\mathcal{O}_{P_3})+h^0(\omega_S)-1=1\\
 h^1(\EE(-k))=0, \forall k\geq4.\\
\end{array}$$
Therefore, the stable rank $2$ bundle $\EE$ has spectrum $\mathcal{X}=\{r_0r_1^3r_2\}=\mathcal{X}_9^{10}.$ Now 
we computer the numbers of generators:\\
$$\begin{array}{l}
 \rho(-3)=h^0(\EE(-3))=1.\\
 \rho(-1)=0\\
 \rho(-k)=0, k\geq4.\\
\end{array}
 $$
(The morphisms 
$H^1(\EE(-i-1))\otimes H^0(\op3(1))\longrightarrow H^1(\EE(-i)) \mbox{ for } i=1 \mbox{ or } i\geq4$ 
are surjective). To computer $\rho(-2)$ we look to the morphism $\delta$ with $m=0$ which 
has rank $1$. This means $\delta$ injective and thus 
$$H^1(\EE(-3))\simeq H^0(\omega_X)\simeq H^0(\mathcal{O}_{P_3}).$$
On the other hand,
$$H^1(\EE(-2))\simeq H^0(\mathcal{O}_{P_2})\oplus H^0(\mathcal{O}_{P_3}(1))\oplus H^0(\omega_S(1)).$$
Therefore,
$$\rho(-2)=\dim\left[\coker\left(H^1(\EE(-3))\otimes H^0(\op3(1))\longrightarrow H^1(\EE(-2))\right)\right]=2.$$
This is sufficient to conclude the existence of the minimal positive Horrocks 
monad $\mathbf{M}$ whose cohomology is the vector bundle $\EE$.

Now if we take $X=P_2\cup P_3$ union of two plane curves $P_2$ and $P_3$ of degree $2$ and 
$3$, respectively, joined by two points $x_0, x_1$ then  $S=\{x_0, x_1\}$ and consequently
$h^0(\omega_{S}(m))=2, m\in\mathbb{Z}$. From Hartshorne--Serre correspondence c.f. Proposition \ref{prop12} there is a rank $2$ 
bundle $\mathcal{F}$ with $c_1(\mathcal{F})=-1, c_2(\mathcal{F})=10$. Furthermore,
\begin{equation}\label{eq-gene}
\begin{array}{l}
 H^1(\mathcal{F}(-3))\simeq H^0(\omega_X)\simeq H^0(\mathcal{O}_{P_3})\oplus H^0(k(x_0))\\
 H^1(\mathcal{F}(-2))\simeq H^0(\mathcal{O}_{P_2})\oplus H^0(\mathcal{O}_{P_3}(1))\oplus H^0(\omega_S(1)),\\
 \end{array}
 \end{equation}
where $k(x_0)$ is the skyscraper sheaf supported at $x_0$ (it could also be the skyscraper sheaf supported at $x_1$).
By repeating the above arguments
$$ h^1(\mathcal{F}(-1))=11,  h^1(\mathcal{F}(-2))=6,  h^1(\mathcal{F}(-3))=2,  h^1(\mathcal{F}(-k))=0, \forall k\geq4,$$
which implies $\X(\mathcal{F})=\{r_0r_1^2r_ 2^2\}=\X_{10}^{10}$. It is easy to see that 
$$\rho(-3)=2, \rho(-1)=0 \mbox{ and } \rho(-k)=0, k\geq4.$$
Revisiting the isomorphisms in \eqref{eq-gene} we compute the generators of minimal degree $2$ obtaining $\rho(-2)=1$. Then, the positive minimal Horrocks monad $\mathbf{M'}$ exists and its cohomology is the vector bundle $\mathcal{F}$.
\end{proof}

\section{The negative minimal Horrocks monads}\label{negative-monads}
The main goal of this section is to prove that the sequence of $c_2=2n$ integers $\mathcal{X}=\{-2^{n-1},-1,0,1^{n-1}\}$ is 
realized as spectrum of a stable rank 2 bundle $\mathcal{E}$ and this one is given necessarily as cohomology of 
a negative minimal Horrocks monad. Let us explicitly display
the negative minimal monad whose cohomology 
is a stable bundle with spectrum $\mathcal{X}$ for each $n\geq4$.

Before this, we thank Nicolae Manolache for alerting us to a mistake made in 
\cite[Proposition 17]{MF2021}. The sequence of integers $\mathcal{X}_6^8=\{-2^3,-1,0,1^3\}$ is 
realized as spectrum of a stable vector bundle on $\mathbb{P}^3$ and the question proposed by 
Hartshorne \cite[(Q2), pg. 806 ]{HR91} continues without fails cases for stable bundles on $\mathbb{P}^3$ with 
odd determinant while for stable bundles with even determinant, Iustin Coanda in 
\cite[Proposition 3.12 and Example 3.13]{coanda2024} shown that this question is not true. The proof of 
the following Proposition is an adaptation of the arguments of \cite[Lemma 2.12]{HR91}. 
\begin{prop}\label{seq-spec}
 The sequence of $c_2=2n$ integers $\mathcal{X}=\{-2^{n-1},-1,0,1^{n-1}\}, n\geq4$, is realized as 
 spectrum of a stable rank 2 bundle on $\mathbb{P}^3$ with $c_1=-1$ and $c_2=2n$. 
\end{prop}
\begin{proof}
Let $X$ be a divisor of type $(1,n-1)$ on a smooth quadric $S$ in $\mathbb{P}^3$; note that $X$ is a rational curve of degree $n$. Let $Y$ be a double structure on $X$ given by an exact sequence of the form 
\begin{equation*}
0\longrightarrow\mathcal{I}_Y\longrightarrow\mathcal{I}_X\longrightarrow\omega_X(3)\longrightarrow0,
\end{equation*}
    and take $\mathcal{E}$ to be given by an extension 
    \begin{equation*}
0\longrightarrow\op3(-1)\longrightarrow\EE\longrightarrow\mathcal{I}_Y\longrightarrow0.
\end{equation*}
 Since $X\neq\emptyset$ we can verify that $\mathcal{E}$ is a stable rank 2 bundle on $\mathbb{P}^3$ and from the above exact 
 sequences, for each integer $l\geq1$,
  $$H^1(\mathcal{E}(-l))\simeq H^1(\mathcal{I}_Y(-l))\simeq  H^0(\omega_X(-l+3)).$$ 
  On the other hand, we have the exact sequence  
  \begin{equation*}
0\longrightarrow\mathcal{O}_S(-2,-2)\longrightarrow\mathcal{O}_S(-1,n-3)\longrightarrow\omega_X\longrightarrow0.
\end{equation*} 
In this manner,
$$h^1(\mathcal{E}(-l))=h^0(\omega_X(-l+3))=\left\{\begin{array}{cc}
    0, & l>2 \\
    n-1, & l=2.
\end{array}
\right.$$
In particular, 
$$n-1=h^1(\mathcal{E}(-2))-h^1(\mathcal{E}(-3))=\#\{k_j\in\mathcal{X}|k_j\geq1\}.$$
Since $c_2(\mathcal{E})=\deg(Y)=2n$ follows that $\mathcal{X}(\mathcal{E})=\mathcal{X}$.
    
\end{proof}

It is important to comment that the arguments in the proof of Proposition \ref{seq-spec} cannot be used to determine 
which minimal Horrocks monads have the bundle $\mathcal{F}$ as cohomology. For this, we use Macaulay2 to 
find the generators of the modules $F_0$ and $F_1$.  
\begin{Lemma}\label{surjective}
If $C$ is a rational and ACM curve, then the evaluation morphism 
$\rho: H^0(\omega_C(2))\otimes H^0(\mathcal{O}_{\mathbb{P}^{3}}(1))\longrightarrow H^0(\omega_C(3))$ is surjective.
\end{Lemma}
\begin{proof}
Since the $C$ is a rational curve of degree $n$ we have the isomorphisms 
$$H^0(\omega_C(2))\simeq H^0(\mathcal{O}_C((2n-2)pts)) \mbox{ and } H^0(\omega_C(3))\simeq H^0(\mathcal{O}_C((3n-3)pts)).$$
On the other hand, the curve $C$ is ACM which implies the surjectivity of the morphism 
$H^0(\mathcal{O}_{\mathbb{P}^{3}}(1))\rightarrow H^0(\mathcal{O}_C(1))$. Therefore, we get the commutative diagram
$$\xymatrix{
 H^0(\omega_C(2))\otimes H^0(\mathcal{O}_{\mathbb{P}^{3}}(1))\ar@{->>}[d]\ar[r]^{\rho}& H^0(\omega_C(3)) \ar[d]^{\simeqd}\\
 H^0(\mathcal{O}_C((2n-2)pts))\otimes H^0(\mathcal{O}_{C}(1))\ar@{->>}[r]& H^0(\mathcal{O}_C((3n-3)pts)). \\
}$$
 This shows that the morphism $\rho$ is surjective.   
\end{proof}
\begin{prop}\label{sobre2}
    If $X$ is a curve of type $(a, b)$ on a nonsingular quadric surface $Q\subset\mathbb{P}^3$, then 
    the morphism $\eta: H^0(\omega_X(k))\otimes H^0(\mathcal{O}_{\mathbb{P}^{3}}(1))\longrightarrow H^0(\omega_X(k+1))$ 
    is surjective for $k\geq 0$.
\end{prop}
\begin{proof}
    We can write $X=\mathcal{C}\cup L_1\cup\cdots\cup L_{b-a}$ where $\mathcal{C}$ is a complete intersection of $Q$ with a surface of degree $a$, and $\{L_1,\dots, L_{b-a\}}$ are mutually disjoint lines in the same family. We get the exact sequence
    $$0\rightarrow\omega_{\mathcal{C}}\rightarrow\omega_X\rightarrow\bigoplus_{i=1}^{b-a}\omega_{L_i}(a)\rightarrow0.$$
    Since the curve $\mathcal{C}$ is ACM follows that $H^1(\omega_C(k))=0$ for $k\geq0$ which implies the surjectivity of the morphism $H^0(\omega_X(k))\rightarrow\bigoplus_{i=1}^{b-a}H^0(\omega_{L_i}(k+a))$. From this subjectivity, we have the commutative diagram
   $$\xymatrix{
 H^0(\omega_C(k))\otimes V\ar@{->>}[d]^{\varphi}\ar[r]&H^0(\omega_X(k))\otimes V\ar[d]^{\eta}\ar@{->>}[r]& \oplus H^0(\omega_{L_i}(k+a))\otimes V \ar@{->>}[d]^{\psi}\\
 H^0(\omega_C(k+1))\ar@{->>}[r]&H^0(\omega_X(k+1))\ar@{->>}[r] &\oplus H^0(\omega_{L_i}(k+1+a)), \\
}$$ 
where $V=H^0(\mathcal{O}_{\mathbb{P}^{3}}(1))$. From Lemma \ref{surjective}, the morphisms $\varphi$ and $\psi$ are 
surjective and so the morphism $\eta$ is surjective.
    \end{proof}

\begin{Theorem}\label{aa1}
Let us consider the sequence of integers $\mathcal{X}_n=\{-2^{n-1},-1,0,1^{n-1}\}, n\geq4$. The sequence $\mathcal{X}_n$ 
is realized as spectrum of a stable bundle $\mathcal{E}$ given as cohomology of the negative minimal monad
\begin{equation*}
\mathcal{A}^{\vee}(-1){\longrightarrow}\mathcal{B}\stackrel{\beta}{\longrightarrow}\mathcal{A},
\end{equation*}
where $\mathcal{A}=\op3(-1)^{\oplus(n-3)}\oplus\op3^(2){\oplus (n-1)}$ and $\mathcal{B}=\op3(-2)^{\oplus (2n-3)}\oplus\op3(1)^{\oplus (2n-3)}$.
\end{Theorem}

\begin{proof}
Let us computer the minimal generators of the Rao module of $\EE$. From Theorem \ref{possibles} follows 
that $\rho(-2)=n-1$ and $\rho(-j)=0,$ for $j=1$ and $j\geq3$. Now let us show that the natural map 
$$\varphi: H^1(\EE(-1))\otimes H^0(\mathcal{O}_{\mathbb{P}^3}(1))\rightarrow H^1(\EE)$$
is surjective and thus $\rho(0)=0$. From the exact sequences \eqref{exact11} and 
\eqref{exact12} the surjectivity of $\varphi$ is equivalent to surjectivity of the map 
$H^1(I_Y(-1)) \otimes H^0(O(1))\rightarrow H^1(I_Y)$ which is surjective if 
and only if the map 
$$\rho: H^0(\omega_X(2)) \otimes H^0(\mathcal{O}_{\mathbb{P}^3}(1))\rightarrow H^0(\omega_X(3))$$
is subjective, where $X$ is a rational curve of type $(1,n-1)$. If we take $k=2$ in Proposition 
\ref{sobre2} then follows that map $\rho$ is subjective. Therefore $\rho(0)=0$ as desired.

On the other hand, with a 
slight adaptation of the results of the \cite[Section 2]{coanda2024} for 
stable rank 2 bundles with odd determinants we get
\begin{equation}\label{neg-terms}
 \rho(i)\leq\max(s(-i-1)-2,0), i\geq1.   
\end{equation}
The above inequality is not true for $i=0$. By applying the inequality \eqref{neg-terms} to the 
spectrum $\mathcal{X}_n$ we obtain 
$$\rho(i)=0, i\geq2 \mbox{ and } \rho(1)\leq n-3.$$
If $\rho(1)=t$ then the Equation \eqref{eq1} implies
$$2n=c_2=(n-1)\cdot2\cdot(2+1)-\displaystyle\sum_{i=1}^{n+t}b_i(b_i+1)\Rightarrow \displaystyle\sum_{i=1}^{n+t}b_i(b_i+1)=4n-6.$$
If $t<n-3$ then at least one of the integers $b_i$ is greater than or equal to $2$ which 
contradicts the minimality of the monad. Therefore, $t=n-3$ and $b_i=1, i=1,\cdots, 2n-3$ is the 
only solution of 
$$\displaystyle\sum_{i=1}^{n+t}b_i(b_i+1)=4n-6.$$

\end{proof}
If we apply the formula \eqref{neg-terms} on each spectrum listed in Table \ref{spectra} then the spectra 
that can be realized as cohomology of a negative minimal monad are: 
$\X_5^{10}=\{r_0^2r_1^3\}, \X_8^{10}=\{r_0r_1^4\}$ and $\X_9^{10}=\{r_0r_1^3r_2\}$. The spectrum $\X_8^{10}$ 
was treated in Theorem \ref{aa1} by taking $n=5$ while the formula \eqref{neg-terms} applied to 
the spectrum $\X_5^{10}$ provides $\rho(1)\leq1$. In this case, the possible negative 
minimal monad is $\boldsymbol{b}=(1,1,1,1,0)$ 
and $\boldsymbol{a}=(2,2,2,-1)$ whose cohomology is not a stable bundle as in next proposition.

\begin{prop}
If the negative minimal Horrocks monad $\mathbf{M}$ exists, with $\boldsymbol{b}=(1,1,1,1,0)$ 
and $\boldsymbol{a}=(2,2,2,-1)$, then its cohomology $\EE$ is not a stable bundle.
\end{prop}
\begin{proof}
Suppose, by contradiction, that $\EE$ is a stable bundle. Denote $B:=\op3(1)^{\oplus4}\oplus\op3\oplus\op3(-1)\oplus\op3(-2)^{\oplus4}$  and we 
have the commutative diagram
\begin{equation}\label{diagram3}
\xymatrix{
 K'\ar@{^{(}->}[d]\ar[r]& \op3(1)^{\oplus4}\oplus\op3 \ar@{^{(}->}[d]\ar[r]^{\psi} & \op3(2)^{\oplus3}\ar@{^{(}->}[d] \\
 K\ar[r]\ar[d]& B \ar[r]^{\beta} \ar@{->>}[d]& \op3(2)^{\oplus3}\oplus\op3(-1)\ar[d].\\
 K''\ar[r]&\op3(-1)\oplus\op3(-2)^{\oplus4}\ar@{->>}[r]^{\phi}&\op3(-1).\\
}
\end{equation}
From diagram \eqref{diagram3} and Euler sequence follows that $K''=\Omega_{\mathbb{P}^3}^1(-1)\oplus\op3(-1)$. On the 
other hand, we also have the commutative diagram
$$\xymatrix{
 & K' \ar@{^{(}->}[d]\ar@{-->}[rd]^{\varphi}& \\
 \op3(-3)^{\oplus3}\oplus\op3\ar[r]& K \ar[r]\ar@{->>}[d]& \EE\\
 &K''&\\
}$$
where $\varphi:K'\longrightarrow K\longrightarrow\EE$ is the composition morphism. Let us consider the possible values for the rank of the morphism $\varphi$.
\begin{itemize}
    \item If the morphism $\varphi: K\hookrightarrow\EE$ is injective, then $h^0(\EE)\geq1$ which is a 
    contradiction because $\EE$ is a stable bundle. 
    \item For $\varphi=0$ we get $K'\hookrightarrow\op3(-3)^{\oplus3}\oplus\op3$  which provides three linearly independent 
    sections in $H^0((K')^\vee(-3))$ that is a contradiction because $H^0((K')^\vee(-3))=0$.
    \item Now we assume $\rank\varphi=1$ that is $\Image(\varphi)=I_C(k)$ for some curve $C$ on 
    $\mathbb{P}^3$ and $k\in\mathbb{Z}$. We have the commutative diagram
\begin{equation}\label{diagram2}
\xymatrix{
 \op3(-k-2)\ar@{^{(}->}[d]^{\tau}\ar[r]& K' \ar@{^{(}->}[d]\ar@{->>}[r] & I_C(k)\ar@{^{(}->}[d] \\
 \op3(-3)^{\oplus3}\oplus\op3\ar[r]& K \ar@{->>}[r] \ar@{^{(}->}[d]& \EE\ar[d].\\
 &K''\ar@{->>}[r]&Q,\\
}\end{equation}
where $Q=\coker\varphi$. Since $\EE$ is stable follows that $h^0(I_C(k))=0$ and so 
$h^0(\op3(-k-2))=1$ which implies $k=-2$. In this manner the morphism $\tau:\op3(-k-2)\longrightarrow\op3(-2)^{\oplus3}\oplus\op3$ becomes $\op3\longrightarrow\op3(-3)^{\oplus3}\oplus\op3$ whose cokernel is $\op3(-3)^{\oplus3}$. Then we have the exact sequence
\begin{equation*}
    \op3(-3)^{\oplus3}\hookrightarrow K''\twoheadrightarrow Q
\end{equation*}
which we dualize to obtain
\begin{equation*}
0\longrightarrow Q^\vee\longrightarrow {K''}^\vee\longrightarrow\op3(3)^{\oplus3}\rightarrow\EXT^1(Q,\op3)\longrightarrow0, 
\end{equation*}
with $\EXT^p(Q,\op3)=0, p>1.$ Furthermore, from the last column of the diagram \eqref{diagram2} we look to the commutative diagram
\begin{equation}
\xymatrix{
 I_C(-2)\ar@{^{(}->}[d]\ar[r]& \EE \ar@{=}[d]\ar[r] & Q\ar[d] \\
 \op3(-2)\ar[r]\ar@{->>}[d]& \EE \ar[r]& I_{C'}(1)\\
 \mathcal{O}_C(-2)&&\\
}
\end{equation}
and from Snake Lemma we obtain 
$$0\longrightarrow\mathcal{O}_C(-2)\longrightarrow Q\longrightarrow I_{C'}(1)\longrightarrow0.$$
Therefore, $Q^\vee\simeq\op3(-1)$. On the other hand, 
$$\EXT^2(Q,\op3)\cong\EXT^2(\mathcal{O}_C(-2), \op3)\cong\omega_C(2),$$ 
which implies $C$ be an empty curve and thus $K'=\op3\oplus\op3(-2)$. As $K''=\Omega_{\mathbb{P}^3}^1(-1)\oplus\op3(-1)$ 
from diagram \eqref{diagram3} we have 
$$B=\op3(2)^{\oplus3}\oplus\op3\oplus\op3(-2)\oplus \op3(-1)^{\oplus2}\oplus\Omega_{\mathbb{P}^3}^1(-1)$$ 
that is a contradiction.
    \end{itemize}
\end{proof}

From the formula \eqref{neg-terms} applied to $\X_9^{10}$ we obtain $\rho(1)\leq1$ and 
$\rho(2)=0$. There are four possibilities of negative minimal monads in this case and we will 
just list them here: $\boldsymbol{a}=(3,2,-1)$ e $\boldsymbol{b}=(2,1,0,0)$, $\boldsymbol{a}=(3,2,-1)$ e $\boldsymbol{b}=(1,1,1,1)$, $\boldsymbol{a}=(3,2,2,-1)$ 
e $\boldsymbol{b}=(2,2,1,0,0)$, $\boldsymbol{a}=(3,2,2,-1)$ e $\boldsymbol{b}=(2,1,1,1,1)$.

\section{Dimension of the Moduli space}\label{section5}
Three irreducible components of the moduli space $\mathcal{B}(-1,10)$ are known: the Hartshorne component 
and two Ein components which we comment in Theorem \ref{dimension}. In this section, we apply the 
semi-continuity of the dimension of the cohomology groups
of coherent sheaves and the formula \eqref{dim formula}, whose proof can be found in \cite[p. 23]{MF2021}, to prove that there is at least one 
more new component in $\mathcal{B}(-1,10)$. The notation in this section is the same employed 
in \cite{MF2021} and we won't repeat it here, let us remember the formula for computer 
the dimension of the set of isomorphism classes 
of stable rank $2$ bundles on $\mathbb{P}^3$ with odd determinant 
which is given as the cohomology of a homotopy-free minimal monad. 
\begin{equation} \label{dim formula}
\dim \mathcal{V}(\boldsymbol{a};\boldsymbol{b}) = \dim H - \dim W - \dim\gl(\mathcal{A}) - \dim G.
\end{equation}
Below, we collect in Table 
\ref{table:cc2=10} each family of minimal homotopy-free monad and its respective dimension. 
	\begin{table}[ht]
		\begin{center}
	\begin{tabular}{ | l | c | c | c |c|c|c|c|} 
\hline
Spectrum & $\boldsymbol{b}$ &  $\boldsymbol{a}$ &$w$&$g$&$s$&$h$&$\dim \mathcal{V}(\boldsymbol{a};\boldsymbol{b})$ \\ 
\hline
$\X^{10}_1$ & $(0,0,0,0,0,0)$ & $(1,1,1,1,1)$ &$200$& $25$& $120$ &$420$&$\mathbf{75}$\\
\hline
$\X^{10}_2$&$(0,0,0,0)$ &$(2,1,1)$& $90$&$13$ &$56$&$232$&$73$\\
\hline
$\X^{10}_3$& $(1,0,0)$& $(2,2)$&$56$& $4$& $65$&$198$&$73$\\
\hline
$\X^{10}_6$& $(1,0)$& $(3)$&$0$&$1$&$40$&$121$&$80$\\
\hline
$\X^{10}_5$& $(1,1,1,1)$& $(2,2,2)$&$168$&$9$&$216$&$468$&$\mathbf{75}$\\
\hline
$\X^{10}_{10}$& $(2,2,1)$& $(3,3)$&$120$&$4$&$271$&$484$&$89$\\
\hline
$\X_{10}^{12}$ & $(4,0)$ & $(5)$&$0$& $1$&$317$&$430$&$112$\\
\hline

	\end{tabular}
\caption{Dimension of the families of homotopy-free minimal Horrocks monads whose 
cohomology is a stable bundle with $c_1=-1$ and $c_2=10$. We set 
$h:=\dim H$, $w:=\dim W$, $g:=\dim \gl(\mathcal{A})$, and $s:=\dim G=\dim S$. The families with 
dimensions equal to the expected one are marked in bold.}
\label{table:cc2=10}
\end{center}
\end{table}
\begin{Theorem}\label{dimension}
The Moduli Scheme $\mathcal{B}(-1,10)$ has at least four components:
\begin{enumerate}
    \item The \textit{Hartshorne component} $M_1$ containing the family 
    $\mathcal{V}(1^5;0^4)$  of expected dimension $75$.  
    \item Two \textit{Ein components} $M_2$ and $M_3$ whose a generic point corresponds 
    to an element of the family $\mathcal{V}(3;1,0)$ and $\mathcal{V}(5;0,4)$, 
    respectively. These components have dimensions $80$ and $112$, respectively.
    \item A new component $M_4$ of dimension at least $89$.
\end{enumerate}
\end{Theorem}
\begin{proof}
We can find in \cite[Section 4]{H78} that $M_1$ is an irreducible component of $\mathcal{B}(-1,10)$ with 
expected dimension $75$. In \cite{Ein88} Ein showed that $M_2$ and $M_3$ are irreducible components of 
dimension $80$ and $112$, respectively. As computed in Table \ref{table:cc2=10}, the family 
$\mathcal{V}(3^2;2^2,1)$ has dimension $89$ and thus it cannot be contained in the components 
$M_1, M_2$. On the other hand, if $\mathcal{E}$ is a generic point of $M_3$ 
then $h^1(\mathcal{E}(-5)=1$ and by inferior semi-continuity 
$h^1(\mathcal{E}(-5))\geq1, \forall \EE\in M_3$ while $h^1(\mathcal{F}(-5))=0$ for a vector 
bundle $\mathcal{F}$ of $\mathcal{V}(2^3;1^4)$ and so $\mathcal{V}(2^3;1^4)$ cannot be contained in $M_3$. 

Now we observe that a bundle $\mathcal{F}$  of $\mathcal{V}(2^3;1^4)$ cannot be contained in the Hartshorne component $M_1$ because a 
generic point of $M_1$ belongs to the family $\mathcal{V}(1^5;0^4)$. Furthermore, by inferior semi-continuity 
$$\begin{array}{c}
  h^1(\mathcal{E}(-3))=1\Rightarrow h^1(\mathcal{E}(-3))\geq1, \forall\EE\in M_2   \\
 h^1(\mathcal{G}(-3))=6\Rightarrow h^1(\mathcal{G}(-3))\geq1, \forall\mathcal{G}\in M_3.      
\end{array} $$
Given $\mathcal{F}\in\mathcal{V}(3^2;2^2,1)$ we see that $h^1(\mathcal{F}(-3))=0$, hence 
$\mathcal{V}(3^2;2^2,1)$ cannot be contained neither in $M_2$ nor in $M_3$.  In summary, the families of vector 
bundles $\mathcal{V}(2^3;1^4)$ and $\mathcal{V}(3^2;2^2,1)$ must be contained in a (new) distinct component of $M_1, M_2, M_3$.
\end{proof}
\begin{Remark}
In Table \ref{table:c2=10} we see that there are two families of Ein bundles with different spectra which 
are cohomology of the Ein monads with $b=(1,0)$ and $a=(3)$. By the inferior semi-continuity, a general point 
of the Ein component $M_2$ is a stable vector bundle with spectrum $\mathcal{X}_6^{10}$.
\end{Remark}
\newpage
\bibliographystyle{amsalpha}

\end{document}